\documentclass[draft]{amsart}

\usepackage{amssymb}
\usepackage[all]{xy}

\newtheorem{thm}{Theorem}[section]
\newtheorem{cor}[thm]{Corollary}
\newtheorem{prop}[thm]{Proposition}
\newtheorem{lem}[thm]{Lemma}
\newtheorem{claim}[thm]{Claim}

\theoremstyle{definition}
\newtheorem{dfn}[thm]{Definition}
\newtheorem{ex}[thm]{Example}

\theoremstyle{remark}
\newtheorem{rem}[thm]{Remark}
\newtheorem{caution}[thm]{Caution}

\numberwithin{equation}{section}

\begin{document}

\title[General heart construction on a triangulated category (II)]{General heart construction on a triangulated category (II): Associated cohomological functor}

\author{Noriyuki ABE}\author{Hiroyuki NAKAOKA}
\address{Graduate School of Mathematical Sciences, The University of Tokyo 
3-8-1 Komaba, Meguro, Tokyo, 153-8914 Japan}

\email[Noriyuki ABE]{abenori@ms.u-tokyo.ac.jp}
\email[Hiroyuki NAKAOKA]{deutsche@ms.u-tokyo.ac.jp}


\begin{abstract}
In the preceding part {\rm (I)} of this paper, we showed that for any torsion pair (i.e., $t$-structure without the shift-closedness) in a triangulated category, there is an associated abelian category, which we call the {\it heart}.
Two extremal cases of torsion pairs are $t$-structures and cluster tilting subcategories.
If the torsion pair comes from a $t$-structure, then its heart is nothing other than the heart of this $t$-structure.
In this case, as is well known, by composing certain adjoint functors, we obtain a cohomological functor from the triangulated category to the heart.
If the torsion pair comes from a cluster tilting subcategory, then its heart coincides with the quotient category of the triangulated category by this subcategory.
In this case, the quotient functor becomes cohomological.
In this paper, we unify these two constructions, to obtain a cohomological functor from the triangulated category, to the heart of any torsion pair.
\end{abstract}

\maketitle

\section{Introduction}
Throughout this paper, we fix a triangulated category $\mathcal{C}$.
For any category $\mathcal{K}$, we write abbreviately $K\in\mathcal{K}$, to indicate that $K$ is an object of $\mathcal{K}$.
For any $K,L\in\mathcal{K}$, let $\mathcal{K}(K,L)$ denote the set of morphisms from $K$ to $L$.
If $\mathcal{M},\mathcal{N}$ are full subcategories of $\mathcal{K}$, then $\mathcal{K}(\mathcal{M},\mathcal{N})=0$ means that $\mathcal{K}(M,N)=0$ for any $M\in\mathcal{M}$ and $N\in\mathcal{N}$.
Similarly, $\mathcal{K}(K,\mathcal{N})=0$ means $\mathcal{K}(K,N)=0$ for any $N\in\mathcal{N}$.

By definition, a {\it torsion pair} is a pair of full additive thick subcategories $(\mathcal{X},\mathcal{Y})$ of $\mathcal{C}$, which satisfies the following \cite{I-Y}.
\begin{enumerate}
\item $\mathcal{C}(\mathcal{X},\mathcal{Y})=0$,
\item For any $C\in\mathcal{C}$, there exists a (not necessarily unique) distinguished triangle
\[ X\rightarrow C\rightarrow Y\rightarrow X[1] \]
satisfying $X\in\mathcal{X}$ and $Y\in\mathcal{Y}$.
\end{enumerate}
Remark that if $(\mathcal{X},\mathcal{Y})$ is a torsion pair, then $\mathcal{X}$ and $\mathcal{Y}$ are mutually orthogonal. Namely, 
an object $C$ in $\mathcal{C}$ satisfies $C\in\mathcal{X}$ (resp.~ $C\in\mathcal{Y}$) if and only if $\mathcal{C}(C,\mathcal{Y})=0$ (resp.~ $\mathcal{C}(\mathcal{X},C)=0$).
For a torsion pair $(\mathcal{X},\mathcal{Y})$, its {\it heart} is defined (see Definition \ref{DefC+-}) by
\[ \underline{\mathcal{H}}=(\mathcal{C}^+\cap\mathcal{C}^-)/(\mathcal{X}[1]\cap\mathcal{Y}), \]
and it was shown in \cite{N} that $\underline{\mathcal{H}}$ becomes an abelian category. In the following two extremal cases, this can be described as follows.
\begin{itemize}
\item If $(\mathcal{X},\mathcal{Y})$ satisfies $\mathcal{Y}\subseteq\mathcal{Y}[1]$, then $(\mathcal{X},\mathcal{Y}[1])$ becomes a $t$-structure \cite{BBD}, and $\underline{\mathcal{H}}$ agrees with the heart $\mathcal{H}=\mathcal{X}\cap\mathcal{Y}[1]$ of this $t$-structure.\\
\item If $(\mathcal{X},\mathcal{Y})$ satisfies $\mathcal{X}[1]=\mathcal{Y}$, then $\mathcal{T}=\mathcal{X}[1]=\mathcal{Y}$ becomes a cluster tilting subcategory \cite{K-R} of $\mathcal{C}$, and $\underline{\mathcal{H}}$ agrees with $\mathcal{C}/\mathcal{T}$.
\end{itemize}

A further observation is that, in each of the above two cases, there is a canonically associated cohomological functor from $\mathcal{C}$, to $\mathcal{H}$ or $\mathcal{C}/\mathcal{T}$, respectively (cf.~ \cite{BBD}, \cite{K-Z}). In this paper, for any torsion pair, we construct a cohomological functor from $\mathcal{C}$ to its heart $\underline{\mathcal{H}}$.

\begin{caution}
Unlike \cite{N}, we use a torsion pair instead of a cotorsion pair.
A cotorsion pair $(\mathcal{U},\mathcal{V})$ in \cite{N} corresponds to a torsion pair $(\mathcal{X},\mathcal{Y})$ in this paper, by $\mathcal{U}=\mathcal{X}[1]$ and $\mathcal{V}=\mathcal{Y}$.
\end{caution}

\section*{Acknowledgement}
The second author wishes to thank Professor Toshiyuki Katsura for his encouragement.
The second author also wishes to thank Professor Osamu Iyama, Professor Kiriko Kato, Professor Bernhard Keller for their useful comments, especially on the terminology.
\section{Preliminaries}

We recall some definitions and results from \cite{N}.

\begin{dfn}\label{DefC+-}
For any torsion pair $(\mathcal{X},\mathcal{Y})$ in $\mathcal{C}$, full subcategories $\mathcal{C}^+$ and $\mathcal{C}^-$ of $\mathcal{C}$ are defined as follows. Put $\mathcal{W}=\mathcal{X}[1]\cap\mathcal{Y}$.
\begin{enumerate}
\item $C\in\mathcal{C}^+$ if and only if there exists a distinguished triangle 
\[ Y\rightarrow W\rightarrow C\rightarrow Y[1] \]
satisfying $W\in\mathcal{W}$ and $Y\in\mathcal{Y}$.
\item $C\in\mathcal{C}^-$ if and only if there exists a distinguished triangle
\[ X\rightarrow C\rightarrow W\rightarrow X[1] \]
satisfying $X\in\mathcal{X}$ and $W\in\mathcal{W}$.
\end{enumerate}
\end{dfn}

\begin{rem}[{\cite[Corollary 3.3 and 3.4]{N}}]
$\ \ $
\begin{enumerate}
\item $C\in\mathcal{C}^+$ if and only if any distinguished triangle 
\[ Y\rightarrow X[1]\rightarrow C\rightarrow Y[1]\quad (X\in\mathcal{X},Y\in\mathcal{Y}) \]
satisfies $X[1]\in\mathcal{W}$.
\item $C\in\mathcal{C}^-$ if and only if any distinguished triangle 
\[ X\rightarrow C\rightarrow Y\rightarrow X[1] \quad (X\in\mathcal{X},Y\in\mathcal{Y}) \]
satisfies $Y\in\mathcal{W}$.
\end{enumerate}
\end{rem}

Let $\underline{\mathcal{C}}=\mathcal{C}/\mathcal{W}$ denote the quotient of $\mathcal{C}$ by $\mathcal{W}$.
Since $\mathcal{X}[1]\supseteq\mathcal{W}$ and $\mathcal{Y}\supseteq\mathcal{W}$, we also have additive full subcategories of $\underline{\mathcal{C}}$
\[ \underline{\mathcal{X}[1]}=(\mathcal{X}[1])/\mathcal{W}\ \ \text{and}\ \ \underline{\mathcal{Y}}=\mathcal{Y}/\mathcal{W}. \]
Put $\mathcal{H}=\mathcal{C}^+\cap\mathcal{C}^-$.
Since $\mathcal{H}\supseteq\mathcal{W}$, we have an additive full subcategory
\[ \underline{\mathcal{H}}=\mathcal{H}/\mathcal{W}\subseteq\underline{\mathcal{C}}, \]
which we call {\it the heart} of $(\mathcal{X},\mathcal{Y})$.
Since $\mathcal{C}^+\cap\mathcal{C}^-=\mathcal{H}\supseteq\mathcal{W}$, we also have additive full subcategories of $\underline{\mathcal{C}}$
\[ \underline{\mathcal{C}}^+=\mathcal{C}^+/\mathcal{W}\quad \text{and}\quad \underline{\mathcal{C}}^-=\mathcal{C}^-/\mathcal{W}. \]
\[
\xy
(-16,0)*+{\underline{\mathcal{H}}}="0";
(0,8)*+{\underline{\mathcal{C}}^+}="2";
(0,-8)*+{\underline{\mathcal{C}}^-}="4";
(16,0)*+{\underline{\mathcal{C}}}="6";
{\ar@{^(->} "0";"2"};
{\ar@{^(->} "0";"4"};
{\ar@{^(->} "2";"6"};
{\ar@{^(->} "4";"6"};
{\ar@{}|\circlearrowright "0";"6"};
\endxy
\]

\begin{prop}\label{PropFromGHC1}
Let $(\mathcal{X},\mathcal{Y})$ be a torsion pair in $\mathcal{C}$.
Then we have the following, for any $C\in\mathcal{C}$, $X\in\mathcal{X}$ and $Y\in\mathcal{Y}$. 
\begin{enumerate}
\item $\underline{\mathcal{C}}(X,C)=\mathcal{C}(X,C)$.
\item $\underline{\mathcal{C}}(C,Y[1])=\mathcal{C}(C,Y[1])$.
\item $\underline{\mathcal{C}}(\underline{\mathcal{X}[1]},\underline{\mathcal{Y}})=0$.
\item $\mathcal{C}^+\supseteq\mathcal{Y}[1]$.
\item $\mathcal{C}^-\supseteq\mathcal{X}$.
\end{enumerate}
\end{prop}
\begin{proof}
The proof can be found in \cite{N}.
\end{proof}

\section{Adjoints and Orthogonality}

\begin{prop}[{\cite[Proposition 3.2]{N}}]
$\ \ $
\begin{enumerate}
\item The inclusion functor $\underline{\mathcal{X}[1]}\hookrightarrow\underline{\mathcal{C}}$ admits an additive right adjoint functor $\sigma_{\mathcal{X}[1]}\colon\underline{\mathcal{C}}\rightarrow\underline{\mathcal{X}[1]}$.
Indeed, for any $C\in\mathcal{C}$ and any distinguished triangle
\[ Y_C\rightarrow X_C[1]\overset{x_C}{\longrightarrow}C\overset{y_C}{\longrightarrow}Y_C[1], \]
$\underline{x_C}\colon X_C[1]\rightarrow C$ gives a coreflection of $C$ along $\underline{\mathcal{X}[1]}\hookrightarrow\underline{\mathcal{C}}$. $($For the definition of $($co-$)$reflections, see \cite{Borceux}.$)$

\item Dually, the inclusion functor $\underline{\mathcal{Y}}\hookrightarrow\underline{\mathcal{C}}$ admits an additive left adjoint functor $\sigma_{\mathcal{Y}}\colon\underline{\mathcal{C}}\rightarrow\underline{\mathcal{Y}}$.
\end{enumerate}
\end{prop}
\begin{proof}
We only show {\rm (1)}.
It suffices to show that for any $X\in\mathcal{X}$ and any $x\in\mathcal{C}(X[1],C)$, there exists $s\in\mathcal{C}(X[1],X_C[1])$ satisfying $\underline{x_C}\circ\underline{s}=\underline{x}$, uniquely in $\underline{\mathcal{C}}(X[1],X_C[1])$.

Existence immediately follows from $\mathcal{C}(X[1],Y_C[1])=0$.
To show the uniqueness, suppose $s\in\mathcal{C}(X[1],X_C[1])$ satisfies $\underline{x_C}\circ\underline{s}=0$.
Then $x_C\circ s$ factors through some $W\in\mathcal{W}$ as in the following diagram.
\[
\xy
(-24,0)*+{Y_C}="0";
(-8,0)*+{X_C[1]}="2";
(8,0)*+{C}="4";
(24,0)*+{Y_C[1]}="6";
(8,12)*+{X[1]}="8";
(24,12)*+{W}="10";
{\ar "0";"2"};
{\ar_<<<<<{x_C} "2";"4"};
{\ar_<<<<<{y_C} "4";"6"};
{\ar_{s} "8";"2"};
{\ar^{w_2} "10";"4"};
{\ar^{w_1} "8";"10"};
{\ar@{}|\circlearrowright "8";"4"};
\endxy
\]
Since $\mathcal{C}(W,Y_C[1])=0$, there exists $w_3\in\mathcal{C}(W,X_C[1])$ such that $x_C\circ w_3=w_2$.
Then $s-w_3\circ w_1$ factors through $Y_C$, which means $\underline{s}=\underline{w_3}\circ\underline{w_1}$ since $\underline{\mathcal{C}}(X[1],Y_C)=0$.
Since $\underline{w_3}\circ\underline{w_1}=0$, we obtain $\underline{s}=0$.
\end{proof}

\begin{cor}\label{CorOrthoC+}
For any $C\in\mathcal{C}$, the following are equivalent.
\begin{enumerate}
\item $C\in\mathcal{C}^+$.
\item $\sigma_{\mathcal{X}[1]}(C)=0$.
\item $\underline{\mathcal{C}}(\underline{\mathcal{X}[1]},C)=0$.
\end{enumerate}
\end{cor}
\begin{proof}
By the definition of $\mathcal{C}^+$, {\rm (1)} and {\rm (2)} are equivalent.
Equivalence of {\rm (2)} and {\rm (3)} follows from the fact that
$\sigma_{\mathcal{X}[1]}(C)=0$ if and only if $\underline{\mathcal{X}[1]}(\underline{\mathcal{X}[1]},\sigma_{\mathcal{X}[1]}(C))=0$ if and only if $\underline{\mathcal{C}}(\underline{\mathcal{X}[1]},C)=0$.
\end{proof}

Dually, we have the following.
\begin{cor}\label{CorOrthoC-}
For any $C\in\mathcal{C}$, the following are equivalent.
\begin{enumerate}
\item $C\in\mathcal{C}^-$.
\item $\sigma_{\mathcal{Y}}(C)=0$.
\item $\underline{\mathcal{C}}(C,\underline{\mathcal{Y}})=0$.
\end{enumerate}
\end{cor}

\begin{prop}[{\cite[Construction 4.2 and Proposition 4.3]{N}}]\label{PropRCR}
$\ \ $
\begin{enumerate}
\item The inclusion functor $i_+\colon\underline{\mathcal{C}}^+\hookrightarrow\underline{\mathcal{C}}$ admits an additive left adjoint functor $\tau^+\colon \underline{\mathcal{C}}\rightarrow\underline{\mathcal{C}}^+$.
We denote the adjunction by $\rho\colon \mathrm{Id}_{\underline{\mathcal{C}}}\rightarrow i_+\circ \tau^+$.

\item The inclusion functor $i_-\colon\underline{\mathcal{C}}^-\hookrightarrow\underline{\mathcal{C}}$ admits an additive right adjoint functor $\tau^-\colon\underline{\mathcal{C}}\rightarrow\underline{\mathcal{C}}^-$.
We denote the adjunction by $\lambda\colon i_-\circ \tau^-\rightarrow\mathrm{Id}_{\underline{\mathcal{C}}}$.
\end{enumerate}
\end{prop}

A closer look at the proof of Proposition 4.3 in \cite{N} leads to the following definition. 

\begin{dfn}\label{DefRCR}
Let $C$ be any object in $\mathcal{C}$.
\begin{enumerate}
\item A {\it reflection triangle} for $C$ is a diagram of the form
\begin{equation}
\xy
(-28,0)*+{X^{\prime}}="0";
(-8,0)*+{C}="2";
(8,0)*+{Z}="4";
(24,0)*+{X^{\prime}[1]}="6";
(-18,-12)*+{X[1]}="8";
(-18,-5)*+{_{\circlearrowright}}="10";
{\ar^{w} "0";"2"};
{\ar^{z} "2";"4"};
{\ar "4";"6"};
{\ar_{x^{\prime}} "0";"8"};
{\ar_{x} "8";"2"};
\endxy
\label{R-triangle}
\end{equation}
satisfying $X,X^{\prime}\in\mathcal{X}$ and $Z\in\mathcal{C}^+$,
where
\[ X^{\prime}\overset{w}{\longrightarrow}C\overset{z}{\longrightarrow}Z\rightarrow X^{\prime}[1] \]
is a distinguished triangle.

\item A {\it coreflection triangle} for $C$ is a diagram of the form
\begin{equation}
\xy
(-24,0)*+{Y^{\prime}}="0";
(-8,0)*+{K}="2";
(8,0)*+{C}="4";
(28,0)*+{Y^{\prime}[1]}="6";
(18,-12)*+{Y}="8";
(18,-5)*+{_{\circlearrowright}}="10";
{\ar^{j} "0";"2"};
{\ar^{k} "2";"4"};
{\ar^{p} "4";"6"};
{\ar_{y} "4";"8"};
{\ar_{y^{\prime}} "8";"6"};
\endxy
\label{CR-triangle}
\end{equation}
satisfying $Y,Y^{\prime}\in\mathcal{Y}$ and $K\in\mathcal{C}^-$, where
\[ Y^{\prime}\overset{j}{\longrightarrow}K\overset{k}{\longrightarrow}C\overset{p}{\longrightarrow} Y^{\prime}[1] \]
is a distinguished triangle.
\end{enumerate}
\end{dfn}

\begin{rem}
$\ \ $
\begin{enumerate}
\item By (the proof of) Proposition 4.3 in \cite{N}, for any reflection triangle $(\ref{R-triangle})$ for $C$, there exists a unique isomorphism $\tau^+(C)\overset{\cong}{\rightarrow}Z$ in $\underline{\mathcal{C}}^+$, compatible with $\underline{z}$ and $\rho_C$. So we abbreviate $Z$ in the above diagram to $\tau^+(C)$. Moreover, it was shown that a restriction triangle always exists for any $C\in\mathcal{C}$.
\item Dually, for any coreflection triangle $(\ref{CR-triangle})$ for $C$, there exists a unique isomorphism $K\overset{\cong}{\rightarrow}\tau^-(C)$ in $\underline{\mathcal{C}}^-$, compatible with $\underline{k}$ and $\lambda_C$. Similarly, we abbreviate $K$ in the above diagram to $\tau^-(C)$. There exists a coreflection triangle for any $C\in\mathcal{C}$.
\end{enumerate}
\end{rem}

\begin{prop}\label{PropOrthoX}
For any $C\in\mathcal{C}$, the following are equivalent.
\begin{enumerate}
\item $C\in\mathcal{X}[1]$.
\item $\underline{\mathcal{C}}(C,\underline{\mathcal{C}}^+)=0$.
\item $\tau^+(C)=0$.
\end{enumerate}
\end{prop}
\begin{proof}
By Corollary \ref{CorOrthoC+}, {\rm (1)} implies {\rm (2)}.
Conversely, assume $\underline{\mathcal{C}}(C,\underline{\mathcal{C}}^+)=0$.
Since $\mathcal{C}^+\supseteq\mathcal{Y}[1]$, we have $\underline{\mathcal{C}}(C,Y[1])=0$ for any $Y\in\mathcal{Y}$.
Since $\underline{\mathcal{C}}(C,Y[1])=\mathcal{C}(C,Y[1])$ by Proposition \ref{PropFromGHC1}, this means $\mathcal{C}(C,\mathcal{Y}[1])=0$, i.e., $C\in\mathcal{X}[1]$.

Equivalence of {\rm (2)} and {\rm (3)} follows from the fact that
$\tau^+(C)=0$ if and only if $\underline{\mathcal{C}}^+(\tau^+(C),\underline{\mathcal{C}}^+)=0$ if and only if $\underline{\mathcal{C}}(C,\underline{\mathcal{C}}^+)=0$.
\end{proof}

Dually, we have the following.
\begin{prop}\label{PropOrthoY}
For any $C\in\mathcal{C}$, the following are equivalent.
\begin{enumerate}
\item $C\in\mathcal{Y}$.
\item $\underline{\mathcal{C}}(\underline{\mathcal{C}}^-,C)=0$.
\item $\tau^-(C)=0$.
\end{enumerate}
\end{prop}

\section{Compatibility of $\tau^+$ and $\tau^-$.}

\begin{lem}\label{LemC-C-}
Let
\[ X\rightarrow A\overset{f}{\longrightarrow} B\rightarrow X[1] \]
be a distinguished triangle such that $X\in\mathcal{X}$.
Then, $A\in\mathcal{C}^-$ if and only if $B\in\mathcal{C}^-$.
In particular, if an object $C\in\mathcal{C}$ belongs to $\mathcal{C}^-$, then we have $\tau^+(C)\in\mathcal{H}$ $($\cite[Lemma 4.6]{N}$)$.
\end{lem}
\begin{proof}
First we assume $A\in\mathcal{C}^-$.
By Corollary \ref{CorOrthoC-}, it suffices to show $\underline{\mathcal{C}}(A,\underline{\mathcal{Y}})=0$.
Let $Y^{\dag}$ be any object in $\mathcal{Y}$, and let $y^{\dag}\in\mathcal{C}(B,Y^{\dag})$ be any morphism.
Since $A\in\mathcal{C}^-$, we have $\underline{\mathcal{C}}(A,Y^{\dag})=0$. Thus $y^{\dag}\circ f$ factors through some $W\in\mathcal{W}$, as in the following diagram.
\[
\xy
(-26,0)*+{X}="0";
(-8,0)*+{A}="2";
(8,0)*+{B}="4";
(24,0)*+{X[1]}="6";
(-8,-12)*+{W}="8";
(8,-12)*+{Y^{\dag}}="10";
{\ar^{} "0";"2"};
{\ar^{f} "2";"4"};
{\ar "4";"6"};
{\ar_{w_1} "2";"8"};
{\ar^{y^{\dag}} "4";"10"};
{\ar_{w_2} "8";"10"};
{\ar@{}|\circlearrowright "2";"10"};
\endxy
\]
Since $\mathcal{C}(X,W)=0$, there exists $w_3\in\mathcal{C}(B,W)$ such that $w_3\circ f=w_1$.
Then $y^{\dag}-w_2\circ w_3$ factors through $X[1]$, and thus $\underline{y^{\dag}}=\underline{w_2}\circ\underline{w_3}$ by $\underline{\mathcal{C}}(X[1],Y^{\dag})=0$\ \ (Proposition \ref{PropFromGHC1}).
Since $\underline{w_2}\circ\underline{w_3}=0$, this means $\underline{y^{\dag}}=0$.

Conversely, assume $B\in\mathcal{C}^-$.
Let $Y^{\dag}$ be any object in $\mathcal{Y}$.
By $\mathcal{C}(X,Y^{\dag})=0$, any morphism $y^{\dag}\in\mathcal{C}(A,Y^{\dag})$ factors through $B$.
Since $\underline{\mathcal{C}}(B,Y^{\dag})=0$, this implies $\underline{y^{\dag}}=0$.
Thus $\underline{\mathcal{C}}(A,\underline{\mathcal{Y}})=0$, which means $A\in\mathcal{C}^-$.
\end{proof}


\begin{prop}\label{Proptau+-}
There exists a natural isomorphism
\[ \eta\colon\tau^+\tau^-\overset{\cong}{\longrightarrow}\tau^-\tau^+, \]
where $\tau^+\tau^-$ and $\tau^-\tau^+$ are the abbreviation of
\begin{eqnarray*}
\tau^+\tau^-&=&(\underline{\mathcal{C}}\overset{\tau^-}{\rightarrow}\underline{\mathcal{C}}^-\overset{i_-}{\hookrightarrow}\underline{\mathcal{C}}\overset{\tau^+}{\rightarrow}\underline{\mathcal{C}}^+\overset{i_+}{\hookrightarrow}\underline{\mathcal{C}}),\\
\tau^-\tau^+&=&(\underline{\mathcal{C}}\overset{\tau^+}{\rightarrow}\underline{\mathcal{C}}^+\overset{i_+}{\hookrightarrow}\underline{\mathcal{C}}\overset{\tau^-}{\rightarrow}\underline{\mathcal{C}}^-\overset{i_-}{\hookrightarrow}\underline{\mathcal{C}}).
\end{eqnarray*}
\end{prop}
\begin{proof}

Since $\tau^-\tau^+(C)\in\mathcal{C}^+$, for any $C\in\mathcal{C}$ there exists a unique morphism $\eta_C\in\mathcal{C}(\tau^+\tau^-(C),\tau^-\tau^+(C))$ such that $\eta_C\circ\rho_{\tau^-(C)}=\tau^-(\rho_C)$. 
Thus it suffices to show $\eta_C$ is an isomorphism, for each $C\in\mathcal{C}$.

Take a reflection triangle and a coreflection triangle for $C$:
\[
\xy
(-28,0)*+{X^{\prime}}="0";
(-8,0)*+{C}="2";
(8,0)*+{\tau^+(C)}="4";
(24,0)*+{X^{\prime}[1]}="6";
(-18,-12)*+{X[1]}="8";
(-18,-5)*+{_{\circlearrowright}}="10";
{\ar^{w} "0";"2"};
{\ar^>>>>>{z} "2";"4"};
{\ar "4";"6"};
{\ar_{x^{\prime}} "0";"8"};
{\ar_{x} "8";"2"};
\endxy
\]
\[
\xy
(-24,0)*+{Y^{\prime}}="0";
(-8,0)*+{\tau^-(C)}="2";
(8,0)*+{C}="4";
(28,0)*+{Y^{\prime}[1]}="6";
(18,-12)*+{Y}="8";
(18,-5)*+{_{\circlearrowright}}="10";
{\ar^>>>>{j} "0";"2"};
{\ar^<<<<{k} "2";"4"};
{\ar^<<<<<<<{q} "4";"6"};
{\ar_{y} "4";"8"};
{\ar_{y^{\prime}} "8";"6"};
\endxy
\]
Since $\underline{z}$ gives a reflection and $Y^{\prime}[1]\in\mathcal{Y}[1]\subseteq\mathcal{C}^+$, there exists $t\in\mathcal{C}(\tau^+(C),Y^{\prime}[1])$ such that $\underline{t}\circ\underline{z}=\underline{q}$ in $\underline{\mathcal{C}}(C,Y^{\prime}[1])$. By $\underline{\mathcal{C}}(C,Y^{\prime}[1])=\mathcal{C}(C,Y^{\prime}[1])$, this means $t\circ z=q$.

If we complete $t$ into a distinguished triangle
\begin{equation}
Y^{\prime}\rightarrow H\overset{h_+}{\longrightarrow}\tau^+(C)\overset{t}{\longrightarrow}Y^{\prime}[1],
\label{H+diag}
\end{equation}
then by the octahedron axiom, we also have a distinguished triangle
\begin{equation}
X^{\prime}\overset{s}{\longrightarrow} \tau^-(C)\overset{h_-}{\longrightarrow}H\rightarrow X^{\prime}[1].
\label{H-diag}
\end{equation}

\begin{claim}\label{ClaimInProof}$\ \ $
\begin{enumerate}
\item $t$ factors through $Y$, and diagram $(\ref{H+diag})$ is a coreflection triangle for $\tau^+(C)$.

\item $s$ factors through $X[1]$, and diagram $(\ref{H-diag})$ is a reflection triangle for $\tau^-(C)$.
\end{enumerate}
\end{claim}
\begin{proof}
Since {\rm (2)} can be shown dually, we only show {\rm (1)}.
By Lemma \ref{LemC-C-}, we have $H\in\mathcal{C}^-$. Thus, it suffices to show that $t$ factors through $Y$.
By $\mathcal{C}(X^{\prime},Y)=0$, there exists $y_0\in\mathcal{C}(\tau^+(C),Y)$ such that $y_0\circ z=y$.
\[
\xy
(-26,0)*+{X^{\prime}}="0";
(-8,0)*+{C}="2";
(8,0)*+{\tau^+(C)}="4";
(24,0)*+{X^{\prime}[1]}="6";
(-8,-12)*+{Y}="8";
(8,-12)*+{Y^{\prime}[1]}="10";
(-3,-4)*+{_{\circlearrowright}}="12";
{\ar "0";"2"};
{\ar^<<<<{z} "2";"4"};
{\ar "4";"6"};
{\ar_{y} "2";"8"};
{\ar^{t} "4";"10"};
{\ar^{y_0} "4";"8"};
{\ar_<<<<<{y^{\prime}} "8";"10"};
\endxy
\]
Then $t-y^{\prime}\circ y_0$ factors through $X^{\prime}[1]$, which means $t=y^{\prime}\circ y_0$, since $\mathcal{C}(X^{\prime}[1],Y^{\prime}[1])=0$.
\end{proof}
%
%
%
%
%
%
%
%
%
Thus we obtain an isomorphism $\tau^+\tau^-(C)\overset{\cong}{\longrightarrow}H\overset{\cong}{\longrightarrow}\tau^-\tau^+(C)$, which is compatible with $\rho_{\tau^-(C)}$ and $\tau^-(\rho_C)$ by construction.
So it must coincide with $\eta_C$, and thus $\eta$ becomes a natural isomorphism.
\end{proof}

\section{Construction of the cohomological functor}

\begin{dfn}\label{DefofH}
By Proposition \ref{Proptau+-}, we have a natural isomorphism of functors
\[ \tau^+\tau^-\circ\underline{(\ )}\cong \tau^-\tau^+\circ\underline{(\ )} \]
where $\underline{(\ )}\colon\mathcal{C}\rightarrow\underline{\mathcal{C}} $
denotes the quotient functor.
Moreover, these functors factors through $\underline{\mathcal{H}}\hookrightarrow\underline{\mathcal{C}}$.

We denote these isomorphic functors abbreviately by
\[ H\colon\mathcal{C}\rightarrow\underline{\mathcal{H}}. \]
\end{dfn}

In the rest of this section, we show $H$ is a cohomological functor (Theorem \ref{ThmCohom}). Remark that $\underline{\mathcal{H}}$ is an abelian category \cite{N}.

\begin{ex}
$\ \ $
\begin{enumerate}
\item If $(\mathcal{X},\mathcal{Y})$ satisfies $\mathcal{Y}\subseteq\mathcal{Y}[1]$ and thus $(\mathcal{X},\mathcal{Y}[1])$ becomes a $t$-structure, then we have
\begin{eqnarray*}
&\underline{\mathcal{C}}^-=\mathcal{C}^-=\mathcal{X},&\\
&\underline{\mathcal{C}}^+=\mathcal{C}^+=\mathcal{Y}[1],&\\
&\underline{\mathcal{H}}=\mathcal{H}=\mathcal{X}\cap\mathcal{Y}[1],&
\end{eqnarray*}
and $\tau^+$ (resp.~ $\tau^-$) is the left (resp.~ right) adjoint of the inclusion functor $\mathcal{X}\hookrightarrow\mathcal{C}$ (resp.~ $\mathcal{Y}[1]\hookrightarrow\mathcal{C}$).
Thus, $H$ agrees with the canonical cohomological functor (cf.~ \cite{BBD})
\[ H=\tau^+\tau^-\cong\tau^-\tau^+\colon\mathcal{C}\rightarrow\mathcal{H}. \]

\item If $(\mathcal{X},\mathcal{Y})=(\mathcal{T},\mathcal{T})$ for a cluster tilting subcategory $\mathcal{T}\subseteq\mathcal{C}$, then we have
\begin{eqnarray*}
&\mathcal{C}^+=\mathcal{C}^-=\mathcal{C},&\\
&\underline{\mathcal{H}}=\mathcal{C}/\mathcal{T},&
\end{eqnarray*}
and $\tau^+=\tau^-=\mathrm{Id}_{\underline{\mathcal{C}}}$. Thus, $H$ agrees with the quotient functor
\[ \underline{(\ \ )}\colon\mathcal{C}\rightarrow\mathcal{C}/\mathcal{T}, \]
which can be shown to be cohomological, by the arguments in \cite{K-Z}.
\end{enumerate}
\end{ex}

\begin{prop}\label{PropVanishH}
Let $(\mathcal{X},\mathcal{Y})$ be a torsion pair in $\mathcal{C}$, and let $H:\mathcal{C}\rightarrow\underline{\mathcal{H}}$ be the additive functor constructed above. Then we have
\[ H(\mathcal{X}[1])=H(\mathcal{Y})=0. \]
\end{prop}
\begin{proof}
This immediately follows from Proposition \ref{PropOrthoX}, \ref{PropOrthoY} and \ref{Proptau+-}.
\end{proof}

\begin{lem}\label{LemEx1}
Let
\[ Y[-1]\overset{f}{\longrightarrow}A\overset{g}{\longrightarrow}B\overset{h}{\longrightarrow}Y \]
be any distinguished triangle satisfying $Y\in\mathcal{Y}$ and $A\in\mathcal{C}^-$.
Then
\[ H(A)\overset{H(g)}{\longrightarrow}H(B)\rightarrow 0 \]
is an exact sequence in $\underline{\mathcal{H}}$.
\end{lem}
\begin{proof}
First we show that we may assume $B\in\mathcal{C}^-$.
Take a coreflection triangle for $B$,
\[
\xy
(-24,0)*+{Y^{\prime}}="0";
(-8,0)*+{\tau^-(B)}="2";
(8,0)*+{B}="4";
(28,0)*+{Y^{\prime}[1]}="6";
(18,-12)*+{Y^{\prime\prime}}="8";
(32,-13)*+{.}="9";
(18,-5)*+{_{\circlearrowright}}="10";
{\ar^>>>>{j} "0";"2"};
{\ar^<<<<{k} "2";"4"};
{\ar^<<<<<<<{p} "4";"6"};
{\ar_{y^{\prime\prime}} "4";"8"};
{\ar_{y^{\prime}} "8";"6"};
\endxy
\]
Since $A$ belongs to $\mathcal{C}^-$, there exists some $g^{\prime}\in\mathcal{C}(A,\tau^-(B))$ such that $g=k\circ g^{\prime}$ by Proposition 4.3 in \cite{N}. If we complete $g^{\prime}$ into a distinguished triangle 
\[ A\overset{g^{\prime}}{\longrightarrow}\tau^-(B)\rightarrow Y_0\rightarrow A[1], \]
then we obtain a morphism of triangles
\[
\xy
(-26,7)*+{Y_0[-1]}="0";
(-8,7)*+{A}="2";
(10,7)*+{\tau^-(B)}="4";
(26,7)*+{Y_0}="6";
(-26,-7)*+{Y[-1]}="10";
(-8,-7)*+{A}="12";
(10,-7)*+{B}="14";
(26,-7)*+{Y.}="16";
{\ar "0";"10"};
{\ar@{=} "2";"12"};
{\ar^{k} "4";"14"};
{\ar "6";"16"};
{\ar "0";"2"};
{\ar "10";"12"};
{\ar^>>>>>{g^{\prime}} "2";"4"};
{\ar_{g} "12";"14"};
{\ar "4";"6"};
{\ar "14";"16"};
{\ar@{}|\circlearrowright "0";"12"};
{\ar@{}|\circlearrowright "2";"14"};
{\ar@{}|\circlearrowright "4";"16"};
\endxy
\]
By the octahedron axiom, we have a distinguished triangle
\[ Y^{\prime}\rightarrow Y_0\rightarrow Y\rightarrow Y^{\prime}[1], \]
and thus $Y_0$ belongs to $\mathcal{Y}$.

Remark that $H(k)$ is an isomorphism. Thus, to show $H(g)$ is an epimorphism, it suffices to show $H(g^{\prime})$ is an epimorphism.

Thus we may assume $A,B\in\mathcal{C}^-$. In this case, for any $Q\in\mathcal{H}$, we have a commutative diagram
\[
\xy
(-20,12)*+{\underline{\mathcal{H}}(H(B),Q)}="10";
(20,12)*+{\underline{\mathcal{H}}(H(A),Q)}="12";
(-20,0)*+{\underline{\mathcal{C}}^+(\tau^+(B),Q)}="20";
(20,0)*+{\underline{\mathcal{C}}^+(\tau^+(A),Q)}="22";
(-20,-12)*+{\underline{\mathcal{C}}(B,Q)}="30";
(20,-12)*+{\underline{\mathcal{C}}(A,Q).}="32";
{\ar^{-\circ H(g)} "10";"12"};
{\ar_{-\circ\tau^+(\underline{g})} "20";"22"};
{\ar_{-\circ\underline{g}} "30";"32"};
{\ar_{\cong} "10";"20"};
{\ar^{\cong} "12";"22"};
{\ar_{\cong} "20";"30"};
{\ar^{\cong} "22";"32"};
{\ar@{}|\circlearrowright "10";"22"};
{\ar@{}|\circlearrowright "20";"32"};
\endxy
\]
By $B\in\mathcal{C}^-$, $Y\in\mathcal{Y}$ and Corollary \ref{CorOrthoC-}, we have $\underline{h}=0$.
By Theorem 2.3 in \cite{K-Z} (cf.~ comments before Fact 2.9 in \cite{N}), $\underline{g}$ is epimorphic in $\underline{\mathcal{C}}$.
In particular $\underline{\mathcal{C}}(B,Q)\overset{-\circ\underline{g}}{\longrightarrow}\underline{\mathcal{C}}(A,Q)$ is injective for any $Q\in\mathcal{H}$. Thus Lemma \ref{LemEx1} follows.
\end{proof}


\begin{lem}\label{LemEx2}
Let
\[ X\overset{f}{\longrightarrow}A\overset{g}{\longrightarrow}B\overset{h}{\longrightarrow}X[1] \]
be any distinguished triangle satisfying $X\in\mathcal{X}$. 
Then the sequence
\[ H(X)\overset{H(f)}{\longrightarrow}H(A)\overset{H(g)}{\longrightarrow}H(B)\rightarrow 0 \]
is exact in $\underline{\mathcal{H}}$.
\end{lem}
\begin{proof}
First we show that we may assume $A,B\in\mathcal{C}^-$. Take a coreflection triangle for $B$,
\[
\xy
(-24,0)*+{Y^{\prime}}="0";
(-8,0)*+{\tau^-(B)}="2";
(8,0)*+{B}="4";
(28,0)*+{Y^{\prime}[1]}="6";
(18,-12)*+{Y}="8";
(32,-13)*+{.}="9";
(18,-5)*+{_{\circlearrowright}}="10";
{\ar^>>>>{j} "0";"2"};
{\ar^<<<<{k} "2";"4"};
{\ar^<<<<<<<{p} "4";"6"};
{\ar_{y} "4";"8"};
{\ar_{y^{\prime}} "8";"6"};
\endxy
\]

If we complete $p\circ g$ into a distinguished triangle
\[ Y^{\prime}\rightarrow L\overset{\ell}{\longrightarrow}A\overset{p\circ g}{\longrightarrow}Y^{\prime}[1], \]
then by the octahedron axiom, we obtain a morphism of distinguished triangles
\[
\xy
(-26,7)*+{X}="0";
(-8,7)*+{L}="2";
(8,7)*+{\tau^-(B)}="4";
(24,7)*+{X[1]}="6";
(-26,-7)*+{X}="10";
(-8,-7)*+{A}="12";
(8,-7)*+{B}="14";
(24,-7)*+{X[1].}="16";
{\ar@{=} "0";"10"};
{\ar_{\ell} "2";"12"};
{\ar^{k} "4";"14"};
{\ar@{=} "6";"16"};
{\ar "0";"2"};
{\ar "10";"12"};
{\ar^<<<<{{}^{\exists}g^{\prime}} "2";"4"};
{\ar_{g} "12";"14"};
{\ar "4";"6"};
{\ar "14";"16"};
{\ar@{}|\circlearrowright "0";"12"};
{\ar@{}|\circlearrowright "2";"14"};
{\ar@{}|\circlearrowright "4";"16"};
\endxy
\]
By Lemma \ref{LemC-C-}, we have $L\in\mathcal{C}^-$.
Thus
\[
\xy
(-24,0)*+{Y^{\prime}}="0";
(-8,0)*+{L}="2";
(8,0)*+{A}="4";
(28,0)*+{Y^{\prime}[1]}="6";
(18,-12)*+{Y}="8";
(18,-5)*+{_{\circlearrowright}}="10";
{\ar^{} "0";"2"};
{\ar^{\ell} "2";"4"};
{\ar^<<<<<<<{p\circ g} "4";"6"};
{\ar_{y\circ g} "4";"8"};
{\ar_{y^{\prime}} "8";"6"};
\endxy
\]
is a coreflection triangle for $A$, and
\[ H(\ell)\colon H(L)\rightarrow H(A) \]
becomes an isomorphism in $\underline{\mathcal{H}}$.
Similarly, $H(k)$ is an isomorphism.
\[
\xy
(-32,7)*+{H(X)}="0";
(-11,7)*+{H(L)}="2";
(11,7)*+{H(\tau^-(B))}="4";
(32,7)*+{H(X[1])}="6";
(-32,-7)*+{H(X)}="10";
(-11,-7)*+{H(A)}="12";
(11,-7)*+{H(B)}="14";
(32,-7)*+{H(X[1])}="16";
{\ar@{=} "0";"10"};
{\ar_{H(\ell)}^{\cong} "2";"12"};
{\ar^{H(k)}_{\cong} "4";"14"};
{\ar@{=} "6";"16"};
{\ar "0";"2"};
{\ar "10";"12"};
{\ar^<<<<<{H(g^{\prime})} "2";"4"};
{\ar_{H(g)} "12";"14"};
{\ar "4";"6"};
{\ar "14";"16"};
{\ar@{}|\circlearrowright "0";"12"};
{\ar@{}|\circlearrowright "2";"14"};
{\ar@{}|\circlearrowright "4";"16"};
\endxy
\]
Thus, replacing $A$ by $L$, and $B$ by $\tau^-(B)$, we may assume $A,B\in\mathcal{C}^-$.
Under this assumption, we show $H(B)$ is the cokernel of $H(X)\overset{H(f)}{\longrightarrow}H(A)$. For any $Q\in\mathcal{H}$, we have a commutative diagram
\[
\xy
(-36,12)*+{\underline{\mathcal{H}}(H(B),Q)}="10";
(0,12)*+{\underline{\mathcal{H}}(H(A),Q)}="12";
(40,12)*+{\underline{\mathcal{H}}(H(X),Q)}="14";
(-36,0)*+{\underline{\mathcal{C}}^+(\tau^+(B),Q)}="20";
(0,0)*+{\underline{\mathcal{C}}^+(\tau^+(A),Q)}="22";
(40,0)*+{\underline{\mathcal{C}}^+(\tau^+(X),Q)}="24";
(-36,-12)*+{\underline{\mathcal{C}}(B,Q)}="30";
(0,-12)*+{\underline{\mathcal{C}}(A,Q)}="32";
(40,-12)*+{\underline{\mathcal{C}}(X,Q).}="34";
{\ar^{-\circ H(g)} "10";"12"};
{\ar^{-\circ H(f)} "12";"14"};
{\ar_{-\circ\tau^+(\underline{g})} "20";"22"};
{\ar_{-\circ\tau^+(\underline{f})} "22";"24"};
{\ar_{-\circ\underline{g}} "30";"32"};
{\ar_{-\circ\underline{f}} "32";"34"};
{\ar_{\cong} "10";"20"};
{\ar^{\cong} "12";"22"};
{\ar^{\cong} "14";"24"};
{\ar_{\cong} "20";"30"};
{\ar^{\cong} "22";"32"};
{\ar^{\cong} "24";"34"};
{\ar@{}|\circlearrowright "10";"22"};
{\ar@{}|\circlearrowright "12";"24"};
{\ar@{}|\circlearrowright "20";"32"};
{\ar@{}|\circlearrowright "22";"34"};
\endxy
\]
So it suffices to show
\[ 0\rightarrow\underline{\mathcal{C}}(B,Q)\overset{-\circ \underline{g}}{\longrightarrow}\underline{\mathcal{C}}(A,Q)\overset{-\circ \underline{f}}{\longrightarrow}\underline{\mathcal{C}}(X,Q) \]
is exact for any $Q\in\mathcal{H}$.

First we show $\underline{\mathcal{C}}(B,Q)\overset{-\circ \underline{g}}{\longrightarrow}\underline{\mathcal{C}}(A,Q)$ is injective.
Let $q\in\mathcal{C}(B,Q)$ be any morphism, and suppose $\underline{q}\circ\underline{g}=0$.
Let
\[ X_A\overset{x_A}{\longrightarrow}A\overset{y_A}{\longrightarrow}Y_A\rightarrow X_A[1] \quad (X_A\in\mathcal{X},Y_A\in\mathcal{W}) \]
be a distinguished triangle.

By $\underline{\mathcal{C}}(X_A,Q)=\mathcal{C}(X_A,Q)$ and $\underline{q}\circ\underline{g}=0$, we have $q\circ g\circ x_A=0$.
Thus there exists $s\in\mathcal{C}(Y_A,Q)$ such that $s\circ y_A=q\circ g$.
Moreover by $\mathcal{C}(X,Y_A)=0$, there exists $t\in\mathcal{C}(B,Y_A)$ such that $t\circ g=y_A$.

\[
\xy
(-30,0)*+{X}="0";
(-10,0)*+{A}="2";
(-10,-14)*+{Y_A}="4";
(-10,12)*+{X_A}="10";
(8,0)*+{B}="12";
(8,-14)*+{Q}="14";
(24,0)*+{X[1]}="16";
(2,-9)*+{}="15";
{\ar^{f} "0";"2"};
{\ar^{g} "2";"12"};
{\ar_{y_A} "2";"4"};
{\ar_{s} "4";"14"};
{\ar^{x_A} "10";"2"};
{\ar^{t} "12";"4"};
{\ar^{q} "12";"14"};
{\ar^{h} "12";"16"};
{\ar@{}|\circlearrowright "2";"15"};
\endxy
\]

Then $q-s\circ t$ factors through $X[1]$.
Since $X[1]\in\mathcal{X}[1]$, $Q\in\mathcal{C}^+$ and $\underline{\mathcal{C}}(\underline{\mathcal{X}[1]},\underline{\mathcal{C}}^+)=0$, we obtain $\underline{q}=\underline{s}\circ\underline{t}$.
Since $\underline{s}=0$ by $Y_A\in\mathcal{W}$, this means $\underline{q}=0$.

To show the exactness of $\underline{\mathcal{C}}(B,Q)\overset{-\circ \underline{g}}{\longrightarrow}\underline{\mathcal{C}}(A,Q)\overset{-\circ \underline{f}}{\longrightarrow}\underline{\mathcal{C}}(X,Q)$, suppose $r\in\mathcal{C}(A,Q)$ satisfies $\underline{r}\circ\underline{f}=0$.
Since $\underline{\mathcal{C}}(X,Q)=\mathcal{C}(X,Q)$, this means $r\circ f=0$, and thus $r$ factors through $g$.
\end{proof}

Dually, we have the following.
\begin{lem}\label{LemEx3}
For any distinguished triangle
\[ Y\overset{f}{\longrightarrow}A\overset{g}{\longrightarrow}B\overset{h}{\longrightarrow}Y[1] \]
satisfying $Y\in\mathcal{Y}$, the sequence
\[ 0\rightarrow H(A)\overset{H(g)}{\longrightarrow}H(B)\overset{H(h)}{\longrightarrow}H(Y[1]) \]
is exact in $\underline{\mathcal{H}}$.
\end{lem}

\begin{thm}\label{ThmCohom}
For any torsion pair $(\mathcal{X},\mathcal{Y})$, the additive functor
\[ H\colon\mathcal{C}\rightarrow\underline{\mathcal{H}} \]
defined in Definition \ref{DefofH} is cohomological.
\end{thm}
\begin{proof}
Let
\[ A\overset{f}{\longrightarrow}B\overset{g}{\longrightarrow}C\overset{h}{\longrightarrow}A[1] \]
be any distinguished triangle.
Take distinguished triangles
\begin{eqnarray*}
X\overset{x}{\longrightarrow}A\overset{y}{\longrightarrow}Y\rightarrow X[1],\\
X\overset{f\circ x}{\longrightarrow}B\overset{b}{\longrightarrow}D\rightarrow X[1]
\end{eqnarray*}
satisfying $X\in\mathcal{X}$ and $Y\in\mathcal{Y}$.
By the octahedron axiom, we obtain another distinguished triangle
\[ Y\overset{p}{\longrightarrow}D\overset{d}{\longrightarrow}C\rightarrow Y[1] \]
and the following diagram.
\[
\xy
(-24,0)*+{X}="0";
(2,0)*+{B}="2";
(14,0)*+{D}="4";
(-10,6)*+{A}="6";
(4,12)*+{Y}="8";
(21,-8)*+{C}="10";
(-10,2)*+{_{\circlearrowright}}="12";
(12.5,-2.5)*+{_{\circlearrowright}}="12";
(3,5)*+{_{\circlearrowright}}="12";
{\ar_{f\circ x} "0";"2"};
{\ar^{b} "2";"4"};
{\ar^{x} "0";"6"};
{\ar^{y} "6";"8"};
{\ar^{f} "6";"2"};
{\ar^{p} "8";"4"};
{\ar^{d} "4";"10"};
{\ar_{g} "2";"10"};
\endxy
\]

By Lemma \ref{LemEx1}, \ref{LemEx2} and \ref{LemEx3}, we have exact sequences
\begin{eqnarray*}
H(X)\overset{H(x)}{\longrightarrow}H(A)\rightarrow 0,\\
H(X)\overset{H(f\circ x)}{\longrightarrow}H(B)\overset{H(b)}{\longrightarrow}H(D)\rightarrow 0,\\
0\rightarrow H(D)\overset{H(d)}{\longrightarrow}H(C).
\end{eqnarray*}
From these, we obtain the exactness of $H(A)\overset{H(f)}{\longrightarrow}H(B)\overset{H(g)}{\longrightarrow}H(C)$.
\end{proof}

We can determine the kernel of $H$ as follows.
\begin{cor}\label{CorKernelofH}
For any $C\in\mathcal{C}$, the following are equivalent.
\begin{enumerate}
\item $H(C)=0$ in $\underline{\mathcal{H}}$.

\item For any $X\in\mathcal{X}$, any morphism $x\in\mathcal{C}(X,C)$ factors through $X_0[1]$, for some $X_0\in\mathcal{X}$.

\item For any $Y\in\mathcal{Y}$, any morphism $v\in\mathcal{C}(C,Y[1])$ factors through some $Y_0\in\mathcal{Y}$.
\end{enumerate}
\end{cor}
\begin{proof}
Since {\rm (3)} is the dual of {\rm (2)}, we only show the equivalence of {\rm (1)} and {\rm (2)}. (A direct verification of the equivalence of {\rm (2)} and {\rm (3)} is not difficult either.)

Assume $H(C)=0$, and let
\[
\xy
(-28,0)*+{X_C^{\prime}}="0";
(-8,0)*+{C}="2";
(8,0)*+{Z}="4";
(24,0)*+{X_C^{\prime}[1]}="6";
(-18,-12)*+{X_C[1]}="8";
(-18,-5)*+{_{\circlearrowright}}="10";
{\ar^{w_C} "0";"2"};
{\ar^{} "2";"4"};
{\ar "4";"6"};
{\ar_{} "0";"8"};
{\ar_{} "8";"2"};
\endxy
\]
be a reflection triangle for $C$.
By Proposition \ref{PropOrthoY}, $H(C)=0$ if and only if $\tau^+(C)\in\underline{\mathcal{Y}}$, which means $Z\in\mathcal{Y}$.
For any $X\in\mathcal{X}$, since $\mathcal{C}(X,Z)=0$, any morphism $x\in\mathcal{C}(X,C)$ factors through $w_C$. Thus $x$ factors through $X_C$.

Conversely, assume {\rm (2)}, and take a distinguished triangle
\[ X\overset{x}{\longrightarrow}C\rightarrow Y\rightarrow X[1]\quad(X\in\mathcal{X},Y\in\mathcal{Y}). \]
Since $x$ factors through $X_0[1]$ for some $X_0\in\mathcal{X}$, we have $H(x)=0$ by Proposition \ref{PropVanishH}.
Thus $H(C)=0$ follows from Proposition \ref{PropVanishH} and Theorem \ref{ThmCohom}.
\end{proof}

\end{document}